\documentclass[11pt, a4paper, twoside]{article}
\usepackage{amsfonts}
\usepackage{mathrsfs}
\usepackage{amsmath, amsthm}
\usepackage{pgf,tikz}
\usepackage{amssymb}
\usepackage{multirow}
\usepackage{array,color}
\usepackage{fancyhdr}
\usepackage[colorlinks, linkcolor=black, anchorcolor=black, citecolor=black]{hyperref}

\numberwithin{equation}{section}

\allowdisplaybreaks

\begin{document}

\fancyhf{}

\fancyhead[OR]{\thepage}

\renewcommand{\headrulewidth}{0pt}
\renewcommand{\thefootnote}{\fnsymbol {footnote}}

\theoremstyle{plain} 
\newtheorem{thm}{\indent\sc Theorem}[section] 
\newtheorem{lem}[thm]{\indent\sc Lemma}
\newtheorem{cor}[thm]{\indent\sc Corollary}
\newtheorem{prop}[thm]{\indent\sc Proposition}
\newtheorem{claim}[thm]{\indent\sc Claim}
\theoremstyle{definition} 
\newtheorem{dfn}[thm]{\indent\sc Definition}
\newtheorem{rem}[thm]{\indent\sc Remark}
\newtheorem{ex}[thm]{\indent\sc Example}
\newtheorem{notation}[thm]{\indent\sc Notation}
\newtheorem{assertion}[thm]{\indent\sc Assertion}
%
%
\numberwithin{equation}{section}
\renewcommand{\proofname}{\indent\sc Proof.} 
\def\C{\mathbb{C}}
\def\R{\mathbb{R}}
\def\Rn{{\mathbb{R}^n}}
\def\M{\mathbb{M}}
\def\N{\mathbb{N}}
\def\Q{{\mathbb{Q}}}
\def\Z{\mathbb{Z}}
\def\F{\mathcal{F}}
\def\L{\mathcal{L}}
\def\S{\mathcal{S}}
\def\supp{\operatorname{supp}}
\def\essi{\operatornamewithlimits{ess\,inf}}
\def\esss{\operatornamewithlimits{ess\,sup}}
\def\dlim{\displaystyle\lim}

\fancyhf{}

\fancyhead[EC]{W. LI, H. Wang, D. Yan}

\fancyhead[EL]{\thepage}

\fancyhead[OC]{Pointwise Convergence of Schr\"{o}dinger Means}

\fancyhead[OR]{\thepage}

\renewcommand{\headrulewidth}{0pt}
\renewcommand{\thefootnote}{\fnsymbol {footnote}}

\title{\textbf{Sharp convergence for sequences of Schr\"{o}dinger means and related generalizations}
\footnotetext {This work is supported by the National Natural Science Foundation of China (No.11871452).}
\footnotetext {{}{2000 \emph{Mathematics Subject
 Classification}: 42B25, 42B37.}}
\footnotetext {{}\emph{Key words and phrases}: Schr\"{o}dinger mean, Pointwise convergence, Maximal functions.} } \setcounter{footnote}{0}
\author{
Wenjuan Li, Huiju Wang, Dunyan Yan}

\date{}
\maketitle

\begin{abstract}
For decreasing sequences $\{t_{n}\}_{n=1}^{\infty}$ converging to zero, we obtain the almost everywhere  convergence results for sequences of Schr\"{o}dinger means $e^{it_{n}\Delta}f$, where $f \in H^{s}(\mathbb{R}^{N}), N\geq 2$. The convergence results are sharp up to the endpoints, and the method can also be applied  to get the convergence results for the fractional Schr\"{o}dinger means  and nonelliptic Schr\"{o}dinger means.
\end{abstract}

\section{Introduction}
 The solution of the Schr\"{o}dinger equation
\begin{equation}\label{Eq1}
\begin{cases}
 i\partial_{t}u(x,t)-\Delta u(x,t) =0 \:\:\:\ x \in \mathbb{R}^{N}, t \in \mathbb{R}^{+},\\
u(x,0)=f \\
\end{cases}
\end{equation}
can be formally written as
\begin{equation}
e^{it\Delta }f(x):= \int_{\mathbb{R}^{N}}{e^{ix \cdot \xi +it|\xi|^{2} } \hat{f} (\xi)d\xi }.
\end{equation}
 The convergence problem of determining the optimal $s$ for which $e^{it \Delta }f$ (called Schr\"{o}dinger means) pointwisely converges to $f$
 whenever $f \in H^{s}(\mathbb{R}^{N})$ as $t$ continuously tends to zero has been studied extensively. The convergence result holds for $s\ge1/4$ when $N=1$ by Carleson \cite{C}, and for $s > \frac{N}{2(N+1)}$  when $N \ge 2$ by Du-Guth-Li \cite{DGL} and Du-Zhang \cite{DZ}. These results are sharp (except the endpoints when $N\ge2$) according to Dahlberg-Kenig \cite{DK} and Bourgain \cite{B}. It is worth to mention that a different counterexample was raised by Luc\`{a}-Rogers \cite{luca} for $N \ge 2$.

In this paper, we consider a related problem: to investigate the convergence properties of $e^{it_{n}\Delta}f$, where $t_{n}$ belongs to some decreasing sequence  $\{t_{n}\}_{n=1}^{\infty}$ converging to zero. One may expect that less regularity on $f$ is enough to ensure convergence in this discrete case. However, when $N=1$ and $t_{n}=1/n$, $n=1,2,\cdots$, Carleson \cite{C} proved that the convergence result holds for $s >1/4$ but fails for $s < \frac{1}{8}$. Indeed, it actually fails for $s < 1/4$ by the counterexample in Dahlberg-Kenig \cite{DK}, a detailed explanation can be found in Section 3 of Lee-Rogers \cite{LR}. Recently, this kind of problem  was further considered by \cite{DS, S1,SS1}.
In particular,  under the assumption that $\{t_{n}\}_{n=1}^{\infty} \in {\ell}^{r,\infty}(\mathbb{N})$, $0<r<\infty$, i.e.,
\begin{equation}
\mathop{sup}_{b>0}b^{r}\sharp \biggl\{n\in \mathbb{N}:t_{n}>b\biggl\} < \infty,
\end{equation}
it follows from \cite{DS} that $e^{it_{n}\Delta}f$ pointwisely converges to $f$ if and only if $s  \ge min\{ \frac{r}{2r+1}, \frac{1}{4} \}$ when $N=1$. But when $N \ge 2$, the  convergence results obtained by \cite{S1,SS1} are far from sharp.
This  open problem will be studied in this article.

We first state the main results on convergence for sequences of Schr\"{o}dinger means, which are sharp up to the endpoints. Then we obtain some generalizations to the fractional Schr\"{o}dinger means $e^{it\Delta^{\frac{a}{2}}}f$ ($1 <a < \infty$) and nonelliptic Schr\"{o}dinger means $e^{it_{n}L}f$, where
\begin{equation}
e^{it_{n}\Delta^{\frac{a}{2}} }f(x):= \int_{\mathbb{R}^{N}}{e^{ix \cdot \xi +it_{n}|\xi|^{a} } \hat{f} (\xi)d\xi },
\end{equation}
and
\begin{equation}
e^{it_{n}L }f(x):= \int_{\mathbb{R}^{N}}{e^{ix \cdot \xi +it_{n}(\xi_{1}^{2}- \xi_{2}^{2}\pm...\pm\xi_{N}^{2}) } \hat{f} (\xi)d\xi }.
\end{equation}

\subsection{Convergence for sequences of Schr\"{o}dinger means}

\begin{thm}\label{main theorem}
Let $N \ge 2$ and $r\in (0,\infty)$.  For any decreasing sequence $\{t_{n}\}_{n=1}^{\infty} \in {\ell}^{r,\infty}(\mathbb{N})$ converging to zero and $\{t_{n}\}_{n=1}^{\infty} \subset (0,1)$,
we have
\begin{equation}\label{Eq1.3}
\lim_{n \rightarrow \infty}e^{it_{n}\Delta}f(x) = f(x) \hspace{0.2cm} a.e.\hspace{0.2cm} x\in \mathbb{R}^N
\end{equation}
whenever $f\in H^s(\mathbb{R}^N)$ and $s> s_{0}=min\{\frac{r}{\frac{N+1}{N}r+1}, \frac{N}{2(N+1)}\}$.
\end{thm}
By standard arguments, it is sufficient to
show the corresponding maximal estimate in $\mathbb{R}^{N}$.

\begin{thm}\label{theorem1.1}
Under the assumptions of Theorem \ref{main theorem}, we have
\begin{equation}\label{Eq1.7+}
\biggl\|\mathop{sup}_{n \in \mathbb{N}} |e^{it_{n}\Delta}f|\biggl\|_{L^{2}(B(0,1))} \leq C\|f\|_{H^s(\mathbb{R}^N)},
\end{equation}
whenever $f\in H^s(\mathbb{R}^N)$ and $s> s_{0}=min\{\frac{r}{\frac{N+1}{N}r+1}, \frac{N}{2(N+1)}\}$, where the constant $C$  does not depend on $f$.
\end{thm}

By translation  invariance in the $x-$direction, $B(0,1)$ in Theorem \ref{theorem1.1} can be replaced by any ball of radius $1$ in $\mathbb{R}^{N}$, which implies Theorem \ref{main theorem}. The convergence result is almost sharp by the Niki\u{s}in-Stein maximal principle and the following fact that
 Theorem \ref{theorem1.1} is sharp up to the endpoints.

\begin{thm}\label{couterexamplethm}
For each $r \in (0,\infty)$,  there exists a sequence $\{t_n\}_{n=1}^{\infty}$ which belongs to ${\ell}^{r,\infty}(\mathbb{N})$,  the corresponding maximal estimate (\ref{Eq1.7+}) fails if $s< s_{0}=min\{\frac{r}{\frac{N+1}{N}r+1}, \frac{N}{2(N+1)}\}$.
\end{thm}

\begin{rem}
One expects that the sparser the time sequences become, the lower the regularity of pointwise convergence requires. Theorem 1.2 and Theorem 1.3 reveal a perhaps surprising phenomenon, namely if $0<r<\frac{N}{N+1}$, there is a gain over the pointwise convergence result from \cite{DGL,DZ,B,luca} when time tends continuously to zero, but not when $r\geq\frac{N}{N+1}$. In fact, such phenomenon has also appeared in one-dimensional case, see \cite{DS}.
\end{rem}
The construction of our counterexample  appeared in Section \ref{counterexample}  is inspired by the work \cite{luca}, which is an alternative proof for Bourgain's counterexample that showed the necessary condition for $\lim_{t\rightarrow 0}e^{it\Delta}f(x)=f(x)$, a.e. $x\in \mathbb{R}^N$.

Next we briefly explain how to prove Theorem \ref{theorem1.1}. Notice that when $\frac{r}{\frac{N+1}{N}r+1} \ge \frac{N}{2(N+1)}$,  Theorem \ref{theorem1.1} follows from the celebrated results by   \cite{DGL} ($N=2$), and  \cite{DZ} ($N \ge 3$). Therefore, we only need to consider the case when $\frac{r}{\frac{N+1}{N}r+1} < \frac{N}{2(N+1)}$, so we always assume that $0<r < \frac{N}{N+1}$ in what follows.

By Littlewood-Paley decomposition and standard argument, we just concentrate on the case when $\textmd{supp}\hat{f} \subset  \{\xi: |\xi| \sim 2^{k}\}$, $k \gg 1$. We consider the maximal function
\[\mathop{sup}_{n \in \mathbb{N}: t_{n} \ge 2^{-\frac{2k}{(N+1)r/N+1}}} |e^{it_{n}\Delta}f|\]
and
\[\mathop{sup}_{n \in \mathbb{N}: t_{n} < 2^{-\frac{2k}{(N+1)r/N+1}}} |e^{it_{n}\Delta}f|,\]
respectively. We deal with the first term by the assumption that the decreasing sequence $\{t_{n}\}_{n=1}^{\infty} \in {\ell}^{r,\infty}(\mathbb{N})$ and Plancherel's theorem. For the second term, since $k< \frac{2k}{\frac{N+1}{N}r+1}< 2k$, the proof can be completed by the following theorem.

\begin{thm}\label{theorem1.4}
Let  $j \in \mathbb{R}$ with $k < j < 2k$. For any $\epsilon >0$, there exists a constant $C_{\epsilon} >0$ such that
\begin{equation}\label{ellEq1.8}
\biggl\|\mathop{sup}_{t \in  (0,2^{-j})}|e^{it\Delta}f|\biggl\|_{L^2(B(0,1))} \leq C_{\epsilon}2^{(2k-j)\frac{N}{2(N+1)} +\epsilon k}\|f\|_{L^{2}(\mathbb{R}^N)},
\end{equation}
for all $f$ with supp $\hat{f} \subset  \{\xi: |\xi| \sim 2^{k}\}$. The constant $C_{\epsilon}$ does not depend on $f$, $j$ and $k$.
\end{thm}

In the case $N=1$, similar result was built in \cite{DS} by $TT^{*}$ argument and stationary phase method. But their method seems not to work well in the higher dimensional case.
 In order to prove Theorem \ref{theorem1.4},  we first observe that (\ref{ellEq1.8}) holds true if spatial variable is restricted to a ball of radius $2^{k-j}$.
Due to references \cite{DGL,DZ}, for any function $g$ with supp $\hat{g} \subset  \{\xi: |\xi| \sim 2^{2k-j}\}$, it holds
\[\biggl\|\mathop{sup}_{t \in (0,2^{-(2k-j)})}|e^{it\Delta}g(x)|\biggl\|_{L^2(B(0,1))} \leq C_{\epsilon}2^{(2k-j)\frac{N}{2(N+1)} + \epsilon k}\|g\|_{L^{2}(\mathbb{R}^N)}.\]
 By scaling, we have
 \begin{equation}\label{Eqglobal}
\biggl\| \mathop{sup}_{t \in (0,2^{-j})}|e^{it\Delta}g|\biggl\|_{_{L^{2}(B(0,2^{k-j}))}}  \leq C_{\epsilon} 2^{(2k-j)\frac{N}{2(N+1)} +\epsilon k} \|g\|_{L^2(\mathbb{R}^N)}
\end{equation}
whenever supp $\hat{g} \subset  \{\xi: |\xi| \sim 2^{k}\}$. Then we obtain the following lemma by translation invariance in the $x$-direction.
\begin{lem}\label{theorem1.5}
When $k < j < 2k$,  for any $\epsilon >0$ and $x_{0}  \in \mathbb{R}^{N}$, there exists a constant $C_{\epsilon} >0$ such that
\begin{equation}\label{Eqglobal}
\biggl\| \mathop{sup}_{t \in (0,2^{-j})}|e^{it\Delta}f|\biggl\|_{_{L^{2}(B(x_{0},2^{k-j}))}}  \leq C_{\epsilon} 2^{(2k-j)\frac{N}{2(N+1)} +\epsilon k} \|f\|_{L^2(\mathbb{R}^N)},
\end{equation}
 whenever supp $\hat{f} \subset  \{\xi: |\xi| \sim 2^{k}\}$. The constant $C_{\epsilon}$ does not depend on $x_{0}$ and $f$.
\end{lem}

Then we can obtain Theorem \ref{theorem1.4} with the help of  Lemma \ref{theorem1.5}, wave packets decomposition  and an  orthogonality  argument. See Section \ref{section 2} below for details.  Moreover, we give the following remark on Theorem \ref{theorem1.4}.

\begin{rem}
We notice that Theorem \ref{theorem1.4} is almost sharp when $j=k$ and $j=2k$. Indeed, when $j = 2k$, Sobolev's embedding implies
 \begin{equation}\label{uniform}
\biggl\|\mathop{sup}_{t \in (0,2^{-2k})}|e^{it\Delta}f(x)|\biggl\|_{L^2(B(0,1))} \leq C\|f\|_{L^{2}(\mathbb{R}^N)}.
\end{equation}
By taking $\hat{f}$ as the  the characteristic function on the set $  \{\xi: |\xi| \sim 2^{k}\}$, it can be observed that the uniform estimate (\ref{uniform}) is optimal.
When $j = k$, it follows from \cite{DGL,DZ} then
 \begin{equation}\label{Eq1.9}
\biggl\|\mathop{sup}_{t \in (0,2^{-k})}|e^{it\Delta}f(x)|\biggl\|_{L^2(B(0,1))} \leq C2^{\frac{N}{2(N+1)}k + \epsilon k}\|f\|_{L^{2}(\mathbb{R}^N)}.
\end{equation}
The above inequality (\ref{Eq1.9}) is  sharp up to the endpoints according to  the counterexample in \cite{B} or \cite{luca}.
 However, the  presence of $2^{\epsilon k}$ on the right hand side of inequality (\ref{ellEq1.8}) leads us to lose the endpoint results in Theorem \ref{theorem1.1}.
\end{rem}

\subsection{Related generalizations}
The method we adopted to prove Theorem \ref{theorem1.1} can be generalized to the fractional case and the nonelliptic case.
Then the corresponding  convergence results follow. We omit most of  details of the proof because they are very similar with that of Theorem \ref{theorem1.1}. Moreover, the sharpness of the result for the nonelliptic case will be proved in Section \ref{necessary condition} below.

Firstly, for the fractional case, we have the following maximal estimate. When $a =2$, it coincides with Theorem \ref{theorem1.1}.

\begin{thm}\label{fractional}
Under the conditions of Theorem \ref{theorem1.1}, for $1 < a < \infty$,
we have
\begin{equation}
\biggl\|\mathop{sup}_{n \in \mathbb{N}} |e^{it_{n}\Delta^{\frac{a}{2}}}f|\biggl\|_{L^{2}(B(0,1))} \leq C\|f\|_{H^s(\mathbb{R}^N)},
\end{equation}
whenever $f\in H^s(\mathbb{R}^N)$ and $s> s_{0}=min\{\frac{a}{2} \cdot \frac{r}{\frac{N+1}{N}r+1}, \frac{N}{2(N+1)}\}$, where the constant $C$  does not depend on $f$.
\end{thm}

Secondly, we introduce the  following maximal estimate for the nonelliptic Schrodinger means. It is sharp up to the endpoints according to the counterexamples stated in Section \ref{necessary condition} below.

\begin{thm}\label{theorem1.6}
Under the conditions of Theorem \ref{theorem1.1},
we have
\begin{equation}
\biggl\|\mathop{sup}_{n \in \mathbb{N}} |e^{it_{n}L}f|\biggl\|_{L^{2}(B(0,1))} \leq C\|f\|_{H^s(\mathbb{R}^N)},
\end{equation}
whenever $f\in H^s(\mathbb{R}^N)$ and $s >s_{0} = min\{\frac{r}{r+1}, \frac{1}{2}\}$, where the constant $C$  does not depend on $f$.
\end{thm}

The proof of  Theorem \ref{theorem1.6} depends heavily on the following theorem.

\begin{thm}\label{theorem1.7}
If supp $\hat{f} \subset  \{\xi : |\xi| \sim \lambda\}$, $\lambda \ge 1$,  then for any small   interval $I$ with
$\lambda^{-2} \le |I| \le \lambda^{-1},$
we have
 \begin{equation}\label{Eq1.8+}
\biggl\|\mathop{sup}_{t \in I}|e^{itL}f(x)|\biggl\|_{L^2(B(0,1))} \leq C\lambda |I|^{\frac{1}{2}}\|f\|_{L^{2}},
\end{equation}
where the constant $C$ does not depend on $f$.
\end{thm}

  Theorem \ref{theorem1.7} follows directly from Sobolev's embedding. Specially,  Theorem \ref{theorem1.7} is sharp when $|I|= \lambda^{-1}$ according to the counterexample in  Rogers-Vargas-Vega \cite{RVV}.
  When $|I|= \lambda^{-2}$, the sharpness can be proved by taking $\tilde{f}$ as the characteristic function over the annulus $  \{\xi : |\xi| \sim \lambda\}$. We point out that the sharpness of Theorem \ref{theorem1.7} enables us to apply the  similar decomposition as Proposition 2.3 in \cite{DS} to get a stronger result than Theorem \ref{theorem1.6} when $r \in (0,1)$.
\begin{thm}\label{theorem2.1}
 If $\{t_{n}\}_{n=1}^{\infty} \in {\ell}^{r(s),\infty}(\mathbb{N})$, $r(s)= \frac{s}{1-s}$.
Then for any $0<s <\frac{1}{2}$, we have
\begin{equation}
\biggl\|\mathop{sup}_{n \in \mathbb{N}} |e^{it_{n}L}f|\biggl\|_{L^{2}(B(0,1))} \leq C\|f\|_{H^s(\mathbb{R}^N)},
\end{equation}
whenever $f\in H^s(\mathbb{R}^N)$, where the constant $C$  does not depend on $f$.
\end{thm}

\begin{rem}
Below, we synthesize our theorems and all results to our best knowledge, and list all almost sharp requirements of regularity on pointwise convergence for different Schr\"{o}dinger-type operators.
\end{rem}
\begin{table*}[htb]
\begin{center}
 \begin{tabular}[l]{>{\centering}p{25mm}|p{25mm}<{\centering}|p{45mm}<{\centering}|p{45mm}<{\centering}}
\hline   \textbf{Operators type} & \textbf{Spatial dimensions} & \textbf{Continuous case $t\rightarrow 0$}& \textbf{Discrete case $t_n\rightarrow 0$}\\
\hline
\multirow{2}{2cm}{Schr\"{o}dinger operator} & $N=1$ & $s\geq \frac{1}{4}$ & $s\geq\min\{\frac{1}{4},\frac{r}{2r+1}\}$\\
\cline {2-4}
& $N\geq 2$ & $s>\frac{N}{2(N+1)}$  & \textcolor{blue}{$s>\min\{\frac{N}{2(N+1)},\frac{r}{\frac{N+1}{N}r+1}\}$}\\
\hline
\multirow{2}{2cm}{Nonelliptic Schr\"{o}dinger} & $N=2$ & $s\geq \frac{1}{2}$  & \textcolor{blue}{$s\geq\min\{\frac{1}{2},\frac{r}{r+1}\}$}\\
\cline {2-4}
 & $N\geq 3$ & $s>\frac{1}{2}$  & \textcolor{blue}{$s>\min\{\frac{1}{2},\frac{r}{r+1}\}$}\\
\hline
\multirow{2}{2cm}{Fractional $a>1$} & $N=1$ & $s\geq \frac{1}{4}$  & $s\geq\min\{\frac{1}{4},\frac{a}{2}\frac{r}{2r+1}\}$\\
\cline {2-4}
 & $N\geq 2$ & $s>\frac{N}{2(N+1)}$ & \textcolor{blue}{$s>\min\{\frac{N}{2(N+1)},\frac{a}{2}\frac{r}{\frac{N+1}{N}r+1}\}$}\\
\hline
\multirow{2}{2cm}{Fractional $0<a<1$} & $N=1$ & $s> \frac{a}{4}$  & $s>\min\{\frac{a}{4},\frac{a}{2}\frac{r}{2r+1}\}$\\
\cline {2-4}
 & $N\geq 2$ & sharp result is open & sharp result is open\\
\hline
\end{tabular}
\end{center}
\end{table*}

In the table above, the results marked in blue come from Theorem \ref{main theorem}, Theorem \ref{fractional}, Theorem \ref{theorem1.6}  and Theorem \ref{theorem2.1} in this paper. For the remaining results, readers can refer to the relevant results of the nonelliptic Schr\"{o}dinger operators in \cite{RVV}; the conclusions about the fractional Schr\"{o}dinger operators when $t$ continuously tends to $0$ can be found in \cite{CK} ($a > 1$) and \cite{W} ($0 < a <1$); other results were introduced at the beginning of the introduction and will not be repeated here.

\textbf{Conventions}: Throughout this article, we shall use the  notation $A\gg B$, which means if there is a sufficiently large constant $G$, which does not depend on the relevant parameters arising in the context in which
the quantities $A$ and $B$ appear, such that $ A\geq GB$. We write $A\sim B$, and mean that $A$ and $B$ are comparable. By
$A\lesssim B$ we mean that $A \le CB $ for some constant $C$ independent of the parameters related to  $A$ and $B$.

\section{Proof of Theorem \ref{theorem1.1} and Theorem \ref{theorem1.4}}\label{section 2}
\textbf{Proof of Theorem \ref{theorem1.1}.}
Let $s_{1}=\frac{r}{\frac{N+1}{N}r+1} + \epsilon$
for some sufficiently small constant $\epsilon >0$.
We decompose $f$ as $f=\sum_{k=0}^{\infty}{f_{k}},$
where $\textmd{supp} \hat{f_{0}} \subset B(0,1)$, $\textmd{supp} \hat{f_{k}} \subset \{\xi: |\xi| \sim 2^{k}\}, k \ge 1$. Then we have
\begin{equation}\label{Eq2.2}
\biggl\|\mathop{sup}_{n \in \mathbb{N}} |e^{it_{n}\Delta}f|\biggl\|_{L^{2}(B(0,1))} \le \sum_{k=0}^{\infty}{\biggl\|\mathop{sup}_{n \in \mathbb{N}} |e^{it_{n}\Delta}f_{k}|\biggl\|_{L^{2}(B(0,1))}}.
\end{equation}

For $k \lesssim 1$ and arbitrary $x \in B(0,1)$, $|e^{it_{n}\Delta}f_{k}(x)| \lesssim \|f_{k}\|_{L^{2}(\mathbb{R}^{N})}$, it is obvious that
\begin{equation}\label{Eq2.3}
\biggl\|\mathop{sup}_{n \in \mathbb{N}} |e^{it_{n}\Delta}f_{k}|\biggl\|_{L^{2}(B(0,1))} \lesssim \|f\|_{H^{s_{1}}(\mathbb{R}^{N})}.
\end{equation}

For each $k \gg 1$, we decompose $\{t_{n}\}_{n=1}^{\infty}$ as
\[A^{1}_{k}:= \biggl\{t_{n}: t_{n} \ge 2^{-\frac{2k}{\frac{N+1}{N}r+1}} \biggl\}\]
and
\[A^{2}_{k}:= \biggl\{t_{n}: t_{n} < 2^{-\frac{2k}{\frac{N+1}{N}r+1}} \biggl\}.\]
Then we have
\begin{align}\label{Eq2.4}
\biggl\|\mathop{sup}_{n \in \mathbb{N}} |e^{it_{n}\Delta}f_{k}|\biggl\|_{L^{2}(B(0,1))}
&\le \biggl\|\mathop{sup}_{n \in \mathbb{N}: t_{n} \in A^{1}_{k}} |e^{it_{n}\Delta}f_{k}|\biggl\|_{L^{2}(B(0,1))}  + \biggl\|\mathop{sup}_{n \in \mathbb{N}: t_{n} \in A^{2}_{k}} |e^{it_{n}\Delta}f_{k}|\biggl\|_{L^{2}(B(0,1))} \nonumber\\
&: =I + II.
\end{align}

We first estimate $I$.
Since  $\{t_{n}\}_{n=1}^{\infty} \in {\ell}^{r,\infty}(\mathbb{N})$, we have
\begin{equation}\label{Eq2.5}
\sharp A^{1}_{k}  \le C 2^{\frac{2rk}{\frac{N+1}{N}r+1}},
\end{equation}
which implies that
\begin{align}\label{Eq2.6}
I &\le \biggl(\sum_{n \in \mathbb{N}: t_{n} \in A^{1}_{k}}{ \biggl\|e^{it_{n}\Delta}f_{k}\biggl\|^{2}_{L^{2}(B(0,1))}}\biggl)^{1/2} \le 2^{\frac{rk}{\frac{N+1}{N}r+1}} \|f_{k}\|_{L^{2}(\mathbb{R}^{N})} \lesssim 2^{-\epsilon k}\|f\|_{H^{s_1}(\mathbb{R}^{N})}.
\end{align}

For $II$, since
\[A^{2}_{k} \subset \biggl(0,  2^{-\frac{2k}{\frac{N+1}{N}r+1}}\biggl).\]
 By previous discussion, we have
$k<\frac{2k}{\frac{N+1}{N}r+1}<2k.$
Then it follows from Theorem \ref{theorem1.4} that,
\begin{equation}\label{Eq2.7}
II \lesssim 2^{(\frac{r}{\frac{N+1}{N}r+1}+\frac{\epsilon}{2})k}\|f_{k}\|_{L^{2}(\mathbb{R}^{N})} \le 2^{-\frac{\epsilon}{2} k}\|f\|_{H^{s_1}(\mathbb{R}^{N})}.
\end{equation}

Inequalities (\ref{Eq2.4}), (\ref{Eq2.6}) and (\ref{Eq2.7}) yield for $k \gg 1$,
\begin{align}\label{Eq2.8}
\biggl\|\mathop{sup}_{n \in \mathbb{N}} |e^{it_{n}\Delta}f_{k}|\biggl\|_{L^{2}(B(0,1))} &\lesssim 2^{-\frac{\epsilon k}{2}}\|f\|_{H^{s_1}(\mathbb{R}^{N})}.
\end{align}

Combining  inequalities (\ref{Eq2.2}), (\ref{Eq2.3}) and (\ref{Eq2.8}),  inequality (\ref{Eq1.7+}) holds true for $s_1$. By the arbitrariness of $\epsilon$, we have finished the proof of Theorem \ref{theorem1.1}.  It remains to prove Theorem \ref{theorem1.4}.

\textbf{ Proof of Theorem \ref{theorem1.4}:} It includes the wave packets decomposition and an orthogonality argument.

$\bullet$  \textbf{Wave packets decomposition. }

We first decompose $e^{it\Delta}f$ on $B(0,1) \times (0,2^{-j})$ in a standard way. For this goal, we decompose the annulus $ \{\xi: |\xi| \sim 2^{k}\}$ into almost disjoint   $2^{j-k}$-cubes $\theta$ with sides parallel to the coordinate axes in $\mathbb{R}^{N}$. Let $2^{k-j}$-cube $\nu$ be dual to $\theta$ and cover $\mathbb{R}^{N}$ by almost disjoint cubes $\nu$. Denote the center of $\theta$ by $c(\theta)$ and the center of $\nu$ by $c(\nu)$. We notice that if $\nu \neq \nu^{\prime}$, then $|c(\nu) - c(\nu^{\prime})| \ge 2^{k-j}$.

Let $\varphi$ be a Schwartz function defined on $\mathbb{R}^{N}$ whose fourier transform is non-negative and supported in a small neighborhood of the origin, and identically equal to $1$ in another smaller interval. Let $\widehat{\varphi_{\theta}}(\xi) = 2^{-\frac{(j-k)N}{2}} \hat{\varphi}(\frac{\xi- c(\theta)}{2^{j-k}})$ and $\widehat{\varphi_{\theta, \nu}}(\xi)= e^{-ic(\nu) \cdot \xi} \widehat{\varphi_{\theta}}(\xi)$. Then $f$ can be decomposed by
\[f = \sum_{\nu} \sum_{\theta} f_{\theta, \nu} = \sum_{\nu} \sum_{\theta}  \langle f, \varphi_{\theta, \nu}\rangle \varphi_{\theta, \nu}, \]
and
\[\|f\|_{L^{2}}^{2} \sim \sum_{\nu} \sum_{\theta}  |\langle f, \varphi_{\theta, \nu}\rangle|^{2}. \]

When $t \in (0,2^{-j})$, integration by parts implies
\[|e^{it\Delta} \varphi_{\theta,\nu}(x)| \le \frac{C_{M} 2^{\frac{(j-k)N}{2}}}{(1+ 2^{j-k}|x-c(\nu) + 2tc(\theta)|)^{M}}. \]
Here $M$ can be sufficiently large. Therefore, $e^{it\Delta} \varphi_{\theta,\nu}(x)$ is essentially supported in a tube
\[T_{\theta, \nu}: = \{(x,t), |x-c(\nu) + 2tc(\theta)| \le 2^{(j-k)(-1+\delta)}, 0 \le t \le 2^{-j} \},\]
where $\delta =\epsilon^3$. The direction of $T_{\theta, \nu}$ is parallel to the vector $(-2c(\theta), 1)$, and the angle between $(-2c(\theta), 1)$ and the $x$-plane is approximately $2^{-k}$.

$\bullet$  \textbf{Orthogonality argument.}

We just give an orthogonality argument under the assumption $j\geq k+\frac{\epsilon k}{N}$. Otherwise, let $j=k+\epsilon_0k$, $0<\epsilon_0<\frac{\epsilon}{N}$, by Lemma \ref{theorem1.5},
\begin{align}
\biggl\|\mathop{sup}_{t
\in  (0,2^{-j})}|e^{it\Delta}f(x)|\biggl\|_{L^2(B(0,1))}
&\leq \biggl(\sum_{m:|x_m|\leq 1}\biggl\|\mathop{sup}_{t
\in  (0,2^{-j})}|e^{it\Delta}f(x)|\biggl\|^2_{L^2(B(x_m, 2^{k-j}))}\biggl)^{1/2} \nonumber\\
&\lesssim 2^{(2k-j)\frac{N}{2(N+1)}+\epsilon k/2+\epsilon_0kN/2}\|f\|_{L^2} \nonumber\\
&\lesssim 2^{(2k-j)\frac{N}{2(N+1)}+\epsilon k}\|f\|_{L^2}.
\end{align}

Now we decompose $B(0,1)$ by $B(0,1) = \cup_{\nu^{\prime}} B(c( \nu^{\prime}), 2^{k-j})$ with $|c(\nu^{\prime})| \lesssim 1$. Then
 \begin{equation}
\biggl\|\mathop{sup}_{t \in  (0,2^{-j})}|e^{it\Delta}f(x)|\biggl\|_{L^2(B(0,1))}^{2} \leq \sum_{ \nu^{\prime}} \biggl\|\mathop{sup}_{t \in  (0,2^{-j})}|e^{it\Delta}f(x)|\biggl\|_{L^2(B(c(\nu^{\prime}),2^{k-j}))}^{2}.
\end{equation}
Fix $c(\nu^{\prime})$, we divide $f$ into two terms
\[f_{1}= \sum_{\theta} \sum_{\nu: |c(\nu)- c(\nu^{\prime})| \le 2^{(j-k)(-1+10 \delta)}} f_{\theta, \nu},\]
and
\[f_{2}= \sum_{\theta} \sum_{\nu: |c(\nu)- c(\nu^{\prime})| > 2^{(j-k)(-1+ 10 \delta)}} f_{\theta, \nu}.\]

For $f_{1}$, by Lemma \ref{theorem1.5} and the $L^{2}$-orthogonality, we  have
\begin{align}
&\sum_{ \nu^{\prime}} \biggl\|\mathop{sup}_{t
\in  (0,2^{-j})}|e^{it\Delta}f_{1}(x)|\biggl\|_{L^2(B(c(\nu^{\prime}),2^{k-j}))}^{2} \nonumber\\
&\le C_{\epsilon} 2^{(2k-j)\frac{N}{N+1} +\epsilon k}  \sum_{ \nu^{\prime}}  \|f_{1}\|_{L^{2}}^{2} \nonumber\\
& \sim  C_{\epsilon} 2^{(2k-j)\frac{N}{N+1} +\epsilon k}  \sum_{ \nu^{\prime}}  \sum_{\theta} \sum_{\nu: |c(\nu)- c(\nu^{\prime})| \le 2^{(j-k)(-1+10 \delta)}} \|f_{\theta, \nu}\|_{L^{2}}^{2} \nonumber\\
& \lesssim C_{\epsilon} 2^{(2k-j)\frac{N}{N+1} +2\epsilon k}   \|f\|_{L^{2}}^{2}.
\end{align}

 We will complete the proof by showing that the contribution from $|e^{it\Delta}f_{2}|$ is negligible when $(x,t)$ belongs to $B(c(\nu^{\prime}), 2^{k-j}) \times (0,2^{-j})$.

Indeed, by Cauchy-Schwartz's inequality and the $L^{2}$-orthogonality, it holds
\begin{align}
|e^{it\Delta}f_{2}| &\le \|f\|_{L^{2}} \biggl( \sum_{\theta} \sum_{\nu: |c(\nu)- c(\nu^{\prime})| > 2^{(j-k)(-1+10 \delta)}} |e^{it\Delta} \varphi_{\theta, \nu}|^{2} \biggl)^{1/2} \nonumber\\
&\le \|f\|_{L^{2}} C_{M} 2^{\frac{(j-k)N}{2}} \biggl( \sum_{\theta} \sum_{\nu: |c(\nu)- c(\nu^{\prime})| > 2^{(j-k)(-1+10 \delta)}} \frac{1}{(1+ 2^{j-k}|x-c(\nu) + 2tc(\theta)|)^{2M}}\biggl)^{1/2}. \nonumber
\end{align}
For each $\theta$, $|x- c(\nu) + 2t c(\theta)| \ge |c(\nu) -c(\nu^{\prime})|/2$, then we have
\begin{align}
& \sum_{\nu: |c(\nu)- c(\nu^{\prime})| > 2^{(j-k)(-1+10 \delta)}} \frac{1}{(1+ 2^{j-k}|x-c(\nu) + 2tc(\theta)|)^{2M}} \nonumber\\
&\le 2^{2M}   \sum_{\substack{l \in \mathbb{N}^{+}\\ l \ge 2^{10\delta\epsilon k /N} }}\sum_{\nu: l2^{k-j} \le |c(\nu)- c(\nu^{\prime})| < (l+1)2^{k-j}}
\frac{1}{(1+ 2^{j-k}|c(\nu)-c(\nu^{\prime})|)^{2M}} \nonumber\\
&\le 2^{2M}   \sum_{\substack{l \in \mathbb{N}^{+}\\ l \ge 2^{10\delta\epsilon k /N} }}\frac{C_{N} l^{N}}{(1+l)^{2M}} \nonumber\\
&\le C_{M,N} 2^{-M\epsilon^4 k}. \nonumber
\end{align}
Notice that the number of $\theta$'s is dominated by $2^{Nk}$. So by choosing $M$ sufficiently large,  for each $(x,t) \in B(c(\nu^{\prime}), 2^{k-j}) \times (0,2^{-j})$, we have
\[ |e^{it\Delta}f_{2}| \le C_{N} 2^{-1000k} \|f\|_{L^{2}}. \]
Then the proof is finished since
\[\sum_{ \nu^{\prime}} \biggl\|\mathop{sup}_{t \in  (0,2^{-j})}|e^{it\Delta}f_{2}(x)|\biggl\|_{L^2(B(c(\nu^{\prime}),2^{k-j}))}^{2} \le C_{N}^{2} 2^{-2000k} \|f\|_{L^{2}}^{2}.\]

\section{A counterexample: Theorem \ref{couterexamplethm}}\label{counterexample}

We notice that the counterexample for $r=\frac{N}{N+1}$ can be also applied to the case when $r > \frac{N}{N+1}$,
since ${\ell}^{N/(N+1),\infty}(\mathbb{N}) \subset {\ell}^{r,\infty}(\mathbb{N})$ and $min\{\frac{r}{\frac{N+1}{N}r+1}, \frac{N}{2(N+1)}\} = \frac{N}{2(N+1)}$ when $r  > \frac{N}{N+1}$. Therefore, next we always assume $r \in (0, \frac{N}{N+1}]$.

Fix $r \in (0, \frac{N}{N+1}]$, we first construct a sequence which belongs to ${\ell}^{r,\infty}(\mathbb{N})$. Put $\beta = \frac{2}{\frac{N+1}{N}r+1}$. Let
$R_{1}=2$ and for each positive integer $n$, $R_{n+1}^{-\beta} \le \frac{1}{2} R_{n}^{-\beta(r+1)}$.
Denote the intervals $I_{n}=[R_{n}^{-\beta(r+1)}, R_{n}^{-\beta})$, $n \in \mathbb{N}^{+}$.
On each $I_{n}$, we get an equidistributed subsequence $ t_{n_{j}}, j =1, 2,..., j_{n} $ such that
\[\{t_{n_{j}}, 1 \le j \le j_{n}\} =: R_{n}^{-\beta(r+1)}\mathbb{Z} \cap I_{n} ,\]
and $t_{n_{j}} - t_{n_{j+1}} = R_{n}^{-\beta(r+1)} $. We claim that the sequence $ t_{n_{j}}, j=1,2,...,j_{n}, n=1,2,...$ belongs to  ${\ell}^{r,\infty}(\mathbb{N})$.

Indeed, according to Lemma 3.2 from \cite{DS}, it suffices to show that
\begin{equation}\label{check sequence}
\mathop{sup}_{b>0}b^{r}\sharp \biggl\{(n,j): b<t_{n_{j}} \le 2b\biggl\} \lesssim 1.
\end{equation}
Notice that we only need to consider $0 <b <1$ because $t_{n_{j}} \in (0,1)$ for each $n$ and $j$. Assume that $(b, 2b] \cap I_{n} \neq \emptyset$ for some $n$, then we have $b < R_{n}^{-\beta}$, $2b \ge R_{n}^{-\beta(r+1)}.$ Therefore,
\[ 2b < 2R_{n}^{-\beta} \le R_{n-1}^{-\beta(r+1)}, \;\:\ b \ge \frac{1}{2} R_{n}^{-\beta(r+1)} \ge R_{n+1}^{-\beta}. \]
This yields $(b, 2b] \cap I_{n^{\prime}} = \emptyset$ for any $n^{\prime} \neq n$,  hence
\[b^{r}\sharp \biggl\{(n,j): b<t_{n_{j}} \le 2b\biggl\} \le b^{r+1}R_{n}^{\beta(r+1)} < 1.\]
Then (\ref{check sequence}) follows by the arbitrariness of $b$.

Our counterexample comes from the following lemma.
\begin{lem}\label{main lemma}
Let $R \gg 1$ and $I=[R^{-\beta(r+1)},  R^{-\beta})$. Assume that the sequence $\{ t_{j} : 1 \le j \le j_{0}\} = R^{-\beta(r+1)}\mathbb{Z} \cap I $ and $t_{j} - t_{j+1} = R^{-\beta(r+1)} $ for each $1 \le j \le j_{0}-1$. Then there exists a function $f$ with supp $\hat{f} \subset B(0,2R)$ such that
\begin{equation}\label{lower bound}
\biggl\|\mathop{sup}_{1 \le j \le j_{0}}|e^{i\frac{t_{j}}{2 \pi}\Delta}f| \biggl \|_{L^{2}(B(0,1))} \gtrsim R^{\frac{1-\beta}{2}} R^{\frac{\beta}{2}} R^{(N-1)(1-\frac{(r+1)\beta}{2})- \epsilon},
\end{equation}
and
\begin{equation}\label{upper bound}
 \|f   \|_{H^{s}(\mathbb{R}^{N})} \lesssim R^{s} R^{\frac{\beta}{4}} R^{\frac{N-1}{2}(1-\frac{(r+1)\beta}{2})}.
\end{equation}
Here $\epsilon >0$ can be sufficiently small.
\end{lem}

Assume that the maximal estimate
\begin{equation}\label{original}
\biggl\|\mathop{sup}_{n}\mathop{sup}_{j} |e^{i\frac{t_{n_{j}}}{2\pi}\Delta}f|\biggl\|_{L^{2}(B(0,1))} \leq C\|f\|_{H^s(\mathbb{R}^N)}
\end{equation}
holds for some $s >0$ and each $f \in H^{s}(\mathbb{R}^{N})$,
then for each $n \in \mathbb{N}^+$,  we have
\begin{equation}\label{reduce}
\biggl\|\mathop{sup}_{j} |e^{i\frac{t_{n_{j}}}{2\pi}\Delta}f|\biggl\|_{L^{2}(B(0,1))} \leq C\|f\|_{H^s(\mathbb{R}^N)}
\end{equation}
whenever $f \in H^{s}(\mathbb{R}^{N})$. Lemma \ref{main lemma} and (\ref{reduce}) yield
\begin{equation}
R_{n}^{\frac{2-\beta}{4}}   R_{n}^{\frac{N-1}{2}(1-\frac{(r+1)\beta}{2})- \epsilon} \le C R_{n}^{s}.
\end{equation}
Then we have $s \ge \frac{r}{\frac{N+1}{N}r +1}$, since $R_{n}$ can be sufficiently large and $\epsilon$ is arbitrarily small. Finally we obtain a sequence $ \frac{t_{n_{j}}}{2\pi}, j=1,2,...,j_{n}, n=1,2,...  \in {\ell}^{r,\infty}(\mathbb{N})$ such that the maximal estimate (\ref{original}) holds only if $s \ge \frac{r}{\frac{N+1}{N}r +1}$.

In the rest of this section, we prove Lemma \ref{main lemma}. Setting
\[\Omega_{1} = \biggl(-\frac{1}{100}R^{\frac{\beta}{2}},  \frac{1}{100}R^{\frac{\beta}{2}} \biggl),\]
\[\Omega_{2}= \biggl\{ \bar{\xi}\in \mathbb{R}^{N-1}:  \bar{\xi} \in 2 \pi R^{\frac{(r+1)\beta}{2}} \mathbb{Z}^{N-1} \cap B(0,R^{1-\epsilon}) \biggl\} + B(0, \frac{1}{1000}),\]
then we define $\hat{f_{1}}(\xi_{1}) = \hat{h}(\xi_{1} + \pi R)$, $\hat{f_{2}}(\bar{\xi}) = \hat{g} (\bar{\xi} + \pi R\theta)$, where $\hat{h} = \chi_{\Omega_{1} }$, $\hat{g} = \chi_{\Omega_{2}}$, and some $\theta \in \mathbb{S}^{N-2}$  (when $N=2$, we denote $\mathbb{S}^0:=(0,1)$) which will be determined later. Define $f$  by $\hat{f} = \hat{f_{1}}\hat{f_{2}}$, it is easy to check that $f$ satisfies (\ref{upper bound}). We are left to prove that inequality (\ref{lower bound}) holds for such $f$. Notice that
\begin{equation}\label{decompose for f}
|e^{i \frac{t_{j}}{2\pi} \Delta}f(x_{1},\bar{x})| = |e^{i\frac{t_{j}}{2\pi} \Delta}f_{1}(x_{1})||e^{i\frac{t_{j}}{2\pi}\Delta}f_{2}(\bar{x})|.
\end{equation}

We first consider $|e^{i\frac{t_{j}}{2\pi}\Delta}f_{1}(x_{1})|$. A change of variables implies
\[ |e^{i\frac{t_{j}}{2\pi} \Delta}f_{1}(x_{1})| = |e^{i\frac{t_{j}}{2\pi}\Delta}h(x_{1}-Rt_{j})| .\]
It is easy to check that $|e^{i\frac{t_{j}}{2\pi}\Delta}h(x_{1})| \gtrsim |\Omega_{1}|$ for each $j$ whenever $|x_{1}| \le R^{-\frac{\beta}{2}}$.  Note that for each $x_{1} \in (0,  R^{1-\beta})$, there exists at least one $t_{j}$ such that
$|x_{1}-  Rt_{j}| \le  R^{1-\beta(r+1)} \le R^{-\frac{\beta}{2}}$ since $\{t_j\}_{j=1}^{j_0}\subset [R^{-\beta(r+1)},R^{-\beta})$ and $t_j-t_{j+1}=R^{-\beta(r+1)}$. Hence we have
\begin{equation}\label{estimate for f1}
 |e^{i\frac{t_{j}}{2\pi}\Delta}f_{1}(x_{1})| \gtrsim |\Omega_{1}|,
 \end{equation}
whenever $x_{1} \in (0, \frac{1}{2}R^{1-\beta})$ and $Rt_{j} \in (x_{1}, x_{1} + R^{-\frac{\beta}{2}})$.

For $|e^{i\frac{t_{j}}{2\pi}\Delta}f_{2}(\bar{x})|$,  we have
\[ |e^{i\frac{t_{j}}{2\pi} \Delta}f_{2}(\bar{x})| = |e^{i\frac{t_{j}}{2\pi} \Delta}g(\bar{x}- Rt_{j} \theta)| .\]
According to Barcel\'{o}-Bennett-Carbery-Ruiz-Vilela \cite{BBCRV}, for each $j$ and $\bar{x} \in U_{0}$,
\begin{equation}\label{estimate for g}
|e^{i\frac{t_{j}}{2\pi} \Delta}g(\bar{x})| \gtrsim |\Omega_{2}|,
\end{equation}
here
\[U_{0}= \biggl\{ \bar{x}\in \mathbb{R}^{N-1}:  \bar{x} \in R^{-\frac{(r+1)\beta}{2}} \mathbb{Z}^{N-1} \cap B(0,2) \biggl\} + B(0, \frac{1}{1000}R^{-1+\epsilon}).\]
We sketch main idea of the proof of inequality (\ref{estimate for g}) for the reader's convenience. Indeed, for each $\bar{\xi} \in \Omega_{2}$, we write $\bar{\xi} =  2 \pi R^{\frac{(r+1)\beta}{2}}l + \bar{\eta}$, $l \in \mathbb{Z}^{N-1} $, $2\pi |l| \le R^{1-\frac{(r+1)\beta}{2}-\epsilon}$, $\bar{\eta} \in B(0, \frac{1}{1000})$.
Then for any $\bar{x}_{m}= R^{-\frac{(r+1)\beta}{2}}m$, $m \in \mathbb{Z}^{N-1}$, $|m| \le 2 R^{\frac{(r+1)\beta}{2}}$,  $t_{j} = R^{-(r+1)\beta}(j_0+1-j)$, $1 \le j \le j_{0}$, we have
\begin{align}
e^{i  \frac{t_{j}}{2\pi} \Delta }g(\bar{x}_{m}) =e^{2\pi i  m \cdot l +  2\pi i (j_0+1-j)|l|^{2}} e^{ i \bar{x}_{m} \cdot \bar{\eta}  + 2i \frac{t_{j}}{2\pi} 2\pi R^{\frac{(r+1)\beta}{2}}l \cdot \bar{\eta} +   i \frac{t_{j}}{2\pi} |\bar{\eta}|^{2} }=e^{ i \bar{x}_{m} \cdot \bar{\eta}  + 2i \frac{t_{j}}{2\pi} 2 \pi R^{\frac{(r+1)\beta}{2}}l \cdot \bar{\eta} +   i \frac{t_{j}}{2\pi} |\bar{\eta}|^{2} }. \nonumber
\end{align}
Noting that $|\bar{x}_{m}| \le 2$, $|t_{j}| \le R^{-\beta}$ and  $|\bar{\eta}| \le \frac{1}{1000}$ imply
\[\biggl | \bar{x}_{m} \cdot \bar{\eta}  + 2 \frac{t_{j}}{2\pi} 2 \pi R^{\frac{(r+1)\beta}{2}}l \cdot \bar{\eta} +    \frac{t_{j}}{2\pi} |\bar{\eta}|^{2}  \biggl| \le \frac{1}{100},\]
then we have
\[|e^{i\frac{t_{j}}{2\pi} \Delta}g(\bar{x}_{m})| \ge \frac{1}{2} |\Omega_{2}|. \]
Moreover, for each $\bar{x} \in U_{0}$, there exits an $\bar{x}_{m}$ such that $|\bar{x} - \bar{x}_{m}| \le \frac{1}{1000} R^{-1 + \epsilon }$, by the mean value theorem and the fact that $|\bar{\xi}| \le 2R^{1-\epsilon}$,
 \[|e^{i\frac{t_{j}}{2\pi} \Delta}g(\bar{x})- e^{i\frac{t_{j}}{2\pi} \Delta}g(\bar{x}_{m})| \le \int_{\mathbb{R}^{N-1}} |\bar{x} -\bar{x}_{m}| |\bar{\xi}| \hat{g}(\bar{\xi}) d\bar{\xi} \le \frac{1}{500} | \Omega_{2} |.\]
Finally we arrive at inequality (\ref{estimate for g})  by the triangle inequality.

Therefore, we have
\begin{equation}\label{estimate for f2}
|e^{i \frac{t_{j}}{2\pi} \Delta}f_{2}(\bar{x})| \gtrsim |\Omega_{2}|,
\end{equation}
if $\bar{x} \in U_{x_{1}} = \bigcup_{j: Rt_{j} \in R^{1-(r+1)\beta}\mathbb{Z} \cap (x_{1}, x_{1} + R^{-\beta/2}) } U_{0} +  Rt_{j} \theta$.
Next we need to select a $\theta \in \mathbb{S}^{N-2}$, such that $|U_{x_{1}}| \gtrsim 1$ for each $x_{1} \in (0, \frac{1}{2} R^{1-\beta})$, which follows if we can prove that there exists a $\theta \in \mathbb{S}^{N-2}$ so that $B(0,1/2)\subset U_{x_1}$ for all $x_{1} \in (0, \frac{1}{2} R^{1-\beta})$. So it remains to prove the claim that there exists a $\theta \in \mathbb{S}^{N-2}$ such that
\[\bigcup_{j: Rt_{j} \in R^{1-\beta(r+1)}\mathbb{Z} \cap (x_{1}, x_{1} + R^{-\beta/2}) } \biggl\{ \bar{x}\in\mathbb{R}^{N-1}:  \bar{x} \in R^{-\frac{(r+1)\beta}{2}} \mathbb{Z}^{N-1} \cap B(0,2) \biggl\} + Rt_{j} \theta\]
is $\frac{1}{1000}R^{-1 + \epsilon}$ dense in the ball $B(0, 1/2)$. In order to apply Lemma 2.1 from Luc\`{a}-Rogers \cite{luca} to get this claim, we first rescale by $R^{\frac{\beta(r+1)}{2}}$, and replace $R^{1+\frac{\beta(r+1)}{2}}t_{j}$ by $s_{j}$, replace $R^{\frac{\beta r}{2}}$ by $R^{\prime}$, recall that $\beta = \frac{2}{\frac{N+1}{N}r+1}$, then we are reduced to show
\[\bigcup_{j: s_{j} \in (R^{\prime})^{1/N} \mathbb{Z} \cap ((R^{\prime})^{(r+1)/r}x_{1}, (R^{\prime})^{(r+1)/r}x_{1} + R^{\prime}) } \biggl\{ \bar{x}:  \bar{x} \in  \mathbb{Z}^{N-1} \cap B(0,2(R^{\prime})^{(r+1)/r}) \biggl\} + s_{j} \theta\]
is $\frac{1}{1000}(R^{\prime})^{-\frac{1}{N} + \frac{(\frac{N+1}{N}r+1)\epsilon}{ r}}$ dense in the ball $B(0, \frac{1}{2} (R^{\prime})^{(r+1)/r})$,
which is equivalent to prove that for any $y\in B(0,\frac{1}{2}(R')^{(r+1)/r})$, there exist
\begin{equation*}
\bar{x}_y\in \mathbb{Z}^{N-1}\cap B(0,2(R')^{(r+1)/r}) \hspace{0.2cm} \textmd{and} \hspace{0.2cm} s_y\in (R^{\prime})^{1/N} \mathbb{Z} \cap ((R^{\prime})^{(r+1)/r}x_{1}, (R^{\prime})^{(r+1)/r}x_{1} + R^{\prime}),
\end{equation*}
such that
\begin{equation*}
|y-\bar{x}_y-s_y\theta)|<\frac{1}{1000}(R^{\prime})^{-\frac{1}{N} + \frac{(\frac{N+1}{N}r+1)\epsilon}{ r}},
\end{equation*}
for a fixed $\theta\in\mathbb{S}^{N-2}$, which is independent of $y$ and $x_1$. This can be implied by the following lemma from Luc\`{a}-Rogers \cite{luca}, here we restate it for reader's convenience.
\begin{lem}[Lemma 2.1, \cite{luca}]\label{lemma1.2}
Let $d\geq 2$, $0<\epsilon,\delta<1$ and $\kappa>\frac{1}{d+1}$. Then, if $\delta<\kappa$ and $R>1$ is sufficiently large, there is $\theta\in \mathbb{S}^{d-1}$ for which, given any $[y]\in \mathbb{T}^d$ and $a\in \mathbb{R}$, there is a $t_y\in R^{\delta}\mathbb{Z}\cap\{a,a+R\}$ such that
\begin{equation*}
|[y]-[t_y\theta]|\leq \epsilon R^{(\kappa-1)/d},
\end{equation*}
where "$[\cdot]$" means taking the quotient $\mathbb{R}^d/\mathbb{Z}^d=\mathbb{T}^d$. Moreover, this remains true with $d=1$, for some $\theta\in (0,1)$.
\end{lem}

We notice that a similar but more detailed proof can be found in Corollary 2.2 of \cite{luca}.

Finally, it follows from inequalities (\ref{decompose for f}), (\ref{estimate for f1}), (\ref{estimate for f2}) that
\[\int_{B(0,1)} \mathop{sup}_{j}|e^{i \frac{t_{j}}{2\pi} \Delta}f(x_{1},\bar{x})|^{2}d\bar{x}dx_{1} \ge \int_{0}^{\frac{R^{1-\beta}}{2}} \int_{U_{x_{1}}} \mathop{sup}_{j}|e^{i \frac{t_{j}}{2\pi} \Delta}f(x_{1},\bar{x})|^{2}d\bar{x}dx_{1} \gtrsim  R^{1-\beta}|\Omega_{1}|^{2} |\Omega_{2}|^{2}, \]
which implies inequality (\ref{lower bound}).


\section{A counterexample for Theorem \ref{theorem1.6}}\label{necessary condition}
For convenience, we first set $N=2$. By changing of variables, the nonelliptic Schr\"{o}dinger operator can be written as
\begin{equation}
e^{it\square}f(x):= \int_{\mathbb{R}^{2}}{e^{ix \cdot \xi +it\xi_{1}\xi_{2}} \hat{f} (\xi)d\xi }.
\end{equation}
For each $r \in (0,1]$,  there exists $\{t_{n}\}_{n=1}^{\infty} \in \ell^{r, \infty}(\mathbb{N})$, such that the maximal estimate
\begin{equation}\label{Eq3.1+}
\biggl\| \mathop{sup}_{n \in \mathbb{N}} |e^{it_{n}\square}f| \biggl\|_{L^{2}(B(0,1))} \le C\|f\|_{H^{s}}
\end{equation}
holds for all $f \in H^{s}(\mathbb{R}^{2})$ only if $s \ge \frac{r}{r+1}$.

Indeed, we choose $\{t_{n}\}_{n=1}^{\infty} \in \ell^{r, \infty}(\mathbb{N})$ but never belongs to $\ell^{r-\epsilon, \infty}(\mathbb{N})$ for any small $\epsilon >0$. Moreover, $t_{n} -t_{n+1}$ is decreasing. According to  Lemma 3.2 in \cite{DS},  we can select $\{b_{j}\}_{j=1}^{\infty}$ and $\{M_{j}\}_{j=1}^{\infty}$ satisfying
$ \lim_{j \rightarrow \infty}b_{j} = 0, \:\ \lim_{j \rightarrow \infty}M_{j} = \infty, $
and
\begin{equation}\label{Eq3.2}
M_{j}b_{j}^{1-r + \epsilon} \le 1,
\end{equation}
such that
\begin{equation}\label{Eq3.2+}
\sharp\biggl\{n: b_{j} < t_{n} \le 2b_{j}\biggl\} \ge M_{j}b_{j}^{-r + \epsilon}.
\end{equation}
By the similar argument as Proposition 3.3 in \cite{DS}, when $t_{n} \le b_{j}$, we have
\begin{equation}\label{Eq3.3}
t_{n}-t_{n+1} \le 2M_{j}^{-1}b_{j}^{r-\epsilon +1}.
\end{equation}

For fixed $j$, choose
$\lambda_{j}=\frac{1}{1000} M_{j}^{\frac{1}{2}}b_{j}^{-\frac{r-\epsilon +1}{2}}$ and
$\widehat{f_{j}}(\xi_{1},\xi_{2}) =\frac{1}{\lambda_{j}} \chi_{[0,\lambda_{j}]\times[{-\lambda_{j}-1,-\lambda_{j}}]}(\xi_{1},\xi_{2}).$
Therefore,
\begin{equation}\label{Eq3.4}
\|f_{j}\|_{H^{\frac{r-\epsilon}{r-\epsilon +1}}}\leq\lambda_{j}^{\frac{r-\epsilon}{r-\epsilon +1}-\frac{1}{2}}.
\end{equation}
Let
$U_{j}=(0,\frac{\lambda_{j}b_{j}}{2}) \times (-\frac{1}{1000}, \frac{1}{1000}).$
Notice that $U_{j} \subset B(0,1)$ due to inequality (\ref{Eq3.2}). Next, we will show that for each $x \in U_{j}$,
\begin{equation}\label{Eq3.5}
\mathop{sup}_{n \in \mathbb{N}} |e^{it_{n}\square}f_{j}| > \frac{1}{2}.
\end{equation}

Changing of variables shows that for each $n \in \mathbb{N}$,
\begin{align}\label{Eq3.6}
|e^{it_{n}\square}f_{j}(x)| =\biggl|\int_{-1}^{0}\int_{0}^{1}{e^{i\lambda_{j}(x_{1}-\lambda_{j}t_{n})\eta_{1}+ix_{2}\eta_{2}+it_{n}\lambda_{j}\eta_{1}\eta_{2}}d\eta_{1}d\eta_{2}}\biggl|.
\end{align}
For each $x \in U_{j}$, there exists a unique $n(x,j)$ such that
\[x_{1} \in (\lambda_{j}t_{n(x,j)+1}, \lambda_{j}t_{n(x,j)}].\]
It is obvious that $t_{n(x,j)+1} \le \frac{b_{j}}{2}$,  then $t_{n(x,j)} \le b_{j}$ due to inequality (\ref{Eq3.2+}) and the assumption that $ t_{n} -t_{n+1}$ is decreasing. Then it follows from inequality (\ref{Eq3.3}) that
\[|\lambda_{j}(x_{1}-\lambda_{j}t_{n(x,j)})\eta_{1}| \le 2\lambda_{j}^{2} M^{-1}_{j}b_{j}^{r-\epsilon+1} \le \frac{1}{1000}.\]
Also,
$| x_{2}\eta_{2}|  \le \frac{1}{1000},$
and by inequality (\ref{Eq3.2}), we have
$|\lambda_{j}t_{n(x,j)}\eta_{1}\eta_{2}|  \le \lambda_{j}b_{j}  \le \frac{1}{1000}.$
Therefore, if we take $n=n(x,j)$ in  (\ref{Eq3.6}), then the phase function will be sufficiently small such
that
$|e^{it_{n(x,j)}\square}f_{j}(x)| > \frac{1}{2}$
for each $x \in U_{j}$, which implies inequality (\ref{Eq3.5}). Then it follows from inequality (\ref{Eq3.4}) and inequality (\ref{Eq3.5}) that
\[\frac{\| \mathop{sup}_{n \in \mathbb{N}} |e^{it_{n}\square}f_j| \|_{L^{2}(B(0,1))}}{\|f_{j}\|_{H^{\frac{r-\epsilon}{r-\epsilon +1}}}} \ge C M_{j}^{\frac{1}{2(r-\epsilon +1)}}.\]
This implies that the maximal estimate (\ref{Eq3.1+}) can not hold when $s \le \frac{r-\epsilon}{r-\epsilon +1}$, hence when
$s < \frac{r}{r+1}$ by the arbitrariness of $\epsilon$.

\begin{rem}
The original idea we adopted to construct the above counterexample  comes from \cite{RVV}. The same idea remains valid in general dimensions. For example, in $\mathbb{R}^{3}$, by changing variables, we can write
\[e^{itL}f(x):= \int_{\mathbb{R}^{3}}{e^{ix \cdot \xi +it(\xi_{1}\xi_{2}\pm \xi_{3}^{2})} \hat{f} (\xi)d\xi }.\]
In order to prove the necessary condition, we only need to take
\[U_{j}=(0,\frac{\lambda_{j}b_{j}}{2}) \times (-\frac{1}{1000}, \frac{1}{1000}) \times (-\frac{1}{1000}, \frac{1}{1000})\]
and
\[\widehat{f_{j}}(\xi_{1},\xi_{2}, \xi_{3}) =\frac{1}{\lambda_{j}} \chi_{[0,\lambda_{j}]\times[{-\lambda_{j}-1,-\lambda_{j}}] \times (0,1)}(\xi_{1},\xi_{2}, \xi_{3}).\]
\end{rem}

\begin{flushleft}
\vspace{0.3cm}\textsc{Wenjuan Li\\School of Mathematics and Statistics\\Northwest Polytechnical University\\710129\\Xi'an, People's Republic of China}

\vspace{0.3cm}\textsc{Huiju Wang\\School of Mathematics and Statistics\\Henan University\\475001\\Kaifeng, People's Republic of China}

\vspace{0.3cm}\textsc{Dunyan Yan\\School of Mathematics Sciences\\University of Chinese Academy of Sciences\\100049\\Beijing, People's Republic of China}

\end{flushleft}

\end{document}